\documentclass[11pt]{amsart}
\usepackage{amsfonts}
\usepackage{amsmath}
\usepackage{amssymb}
\usepackage{amscd}

\input{tcilatex}
\theoremstyle{definition}
\theoremstyle{remark}
\numberwithin{equation}{section}

\begin{document}
\author{D. S. Lubinsky}
\address{School of Mathematics\\
Georgia Institute of Technology\\
Atlanta, GA 30332-0160\\
USA.\\
lubinsky@math.gatech.edu}
\title[Universality Limits]{A New Approach to Universality Limits involving
Orthogonal Polynomials}
\date{January 10, 2007}

\begin{abstract}
We show how localization and smoothing techniques can be used to establish
universality in the bulk of the spectrum for a fixed positive measure $\mu $
on $\left[ -1,1\right] $. Assume that $\mu $ is a regular measure, and is
absolutely continuous in an open interval containing some closed subinterval 
$J$ of $\left( -1,1\right) $. Assume that in $J$, the absolutely continuous
component $\mu ^{\prime }$ is positive and continuous. Then universality in $%
J$ for $\mu $ follows from universality for the classical Legendre weight.
We also establish universality in an $L_{p}$ sense under weaker assumptions
on $\mu .$
\end{abstract}
\maketitle

\section{Introduction and Results\protect\footnote{%
Research supported by NSF grant DMS0400446 and US-Israel BSF grant 2004353}}

Let $\mu $ be a finite positive Borel measure on $\left( -1,1\right) $. Then
we may define orthonormal polynomials%
\begin{equation*}
p_{n}\left( x\right) =\gamma _{n}x^{n}+...,\gamma _{n}>0,
\end{equation*}%
$n=0,1,2,...$ satisfying the orthonormality conditions%
\begin{equation*}
\int_{-1}^{1}p_{n}p_{m}d\mu =\delta _{mn}.
\end{equation*}%
These orthonormal polynomials satisfy a recurrence relation of the form%
\begin{equation}
xp_{n}\left( x\right) =a_{n+1}p_{n+1}\left( x\right) +b_{n}p_{n}\left(
x\right) +a_{n}p_{n-1}\left( x\right) ,
\end{equation}%
where%
\begin{equation*}
a_{n}=\frac{\gamma _{n-1}}{\gamma _{n}}>0\text{ and }b_{n}\in \mathbb{R}%
\text{, }n\geq 1,
\end{equation*}%
and we use the convention $p_{-1}=0$. Throughout we use 
\begin{equation*}
w=\frac{d\mu }{dx}
\end{equation*}%
to denote the absolutely continuous part of $\mu $. A classic result of E.A.
Rakhmanov \cite{Simon2005} asserts that if $w>0$ a.e. in $\left[ -1,1\right] 
$, then $\mu $ belongs to the Nevai-Blumenthal class $\mathcal{M}$, that is 
\begin{equation}
\lim_{n\rightarrow \infty }a_{n}=\frac{1}{2}\text{ and }\lim_{n\rightarrow
\infty }b_{n}=0.
\end{equation}%
We note that there are pure jump and pure singularly continuous measures in $%
\mathcal{M}$, despite the fact that one tends to associate it with weights
that are positive a.e. A class of measures that contains $\mathcal{M}$ is
the class of \textit{regular measures} on $\left[ -1,1\right] $ \cite%
{StahlTotik1992}, defined by the condition 
\begin{equation*}
\lim_{n\rightarrow \infty }\gamma _{n}^{1/n}=\frac{1}{2}.
\end{equation*}

Orthogonal polynomials play an important role in random matrix theory \cite%
{Deift1999}, \cite{Mehta1991}. One of the key limits there involves the
reproducing kernel%
\begin{equation}
K_{n}\left( x,y\right) =\sum_{k=0}^{n-1}p_{k}\left( x\right) p_{k}\left(
y\right) .
\end{equation}%
Because of the Christoffel-Darboux formula, it may also be expressed as%
\begin{equation}
K_{n}\left( x,y\right) =a_{n}\frac{p_{n}\left( x\right) p_{n-1}\left(
y\right) -p_{n-1}\left( x\right) p_{n}\left( y\right) }{x-y}.
\end{equation}%
Define the normalized kernel 
\begin{equation}
\widetilde{K}_{n}\left( x,y\right) =w\left( x\right) ^{1/2}w\left( y\right)
^{1/2}K_{n}\left( x,y\right) .
\end{equation}%
The simplest case of the universality law is the limit%
\begin{equation}
\lim_{n\rightarrow \infty }\frac{\widetilde{K}_{n}\left( x+\frac{a}{%
\widetilde{K}_{n}\left( x,x\right) },x+\frac{b}{\widetilde{K}_{n}\left(
x,x\right) }\right) }{\widetilde{K}_{n}\left( x,x\right) }=\frac{\sin \pi
\left( a-b\right) }{\pi \left( a-b\right) }.
\end{equation}%
Typically this holds uniformly for $x$ in a compact subinterval of $\left(
-1,1\right) $ and $a,b$ in compact subsets of the real line. Of course, when 
$a=b$, we interpret $\frac{\sin \pi \left( a-b\right) }{\pi \left(
a-b\right) }$ as $1$. We cannot hope to survey the vast body of results on
universality limits here - the reader may consult \cite{Baiketal2006}, \cite%
{Deift1999}, \cite{Deiftetal1999}, \cite{Mehta1991} and the forthcoming
proceedings of the conference devoted to the 60th birthday of Percy Deift.

Our goal here is to present what we believe is a new approach, based on
localization and smoothing. Our main result is:\newline
\newline
\textbf{Theorem 1.1} \newline
\textit{Let }$\mu $\textit{\ be a finite positive Borel measure on }$\left(
-1,1\right) $\textit{\ that is regular. Let} $I$\textit{\ be a closed
subinterval of }$\left( -1,1\right) $\textit{\ such that }$\mu $\textit{\ is
absolutely continuous in an open interval containing }$I.$\textit{\ Assume
that }$w$\textit{\ is positive and continuous in }$I$\textit{. Then
uniformly for }$x\in I$\textit{\ and }$a,b$\textit{\ in compact subsets of
the real line, we have}%
\begin{equation}
\lim_{n\rightarrow \infty }\frac{\widetilde{K}_{n}\left( x+\frac{a}{%
\widetilde{K}_{n}\left( x,x\right) },x+\frac{b}{\widetilde{K}_{n}\left(
x,x\right) }\right) }{\widetilde{K}_{n}\left( x,x\right) }=\frac{\sin \pi
\left( a-b\right) }{\pi \left( a-b\right) }.
\end{equation}

Note that we allow the case where $I$ consists of just a single point.%
\newline
\newline
\textbf{Corollary 1.2}\newline
\textit{Let }$m\geq 1$\textit{\ and }%
\begin{equation*}
R_{m}\left( y_{1},y_{2},...,y_{m}\right) =\det \left( \widetilde{K}%
_{n}\left( y_{i},y_{j}\right) \right) _{i,j=1}^{m}
\end{equation*}%
\textit{denote the }$m-$\textit{point correlation function. Uniformly for }$x%
\mathit{\in }I$\textit{, and for given }$\left\{ \xi _{j}\right\} _{j=1}^{m}$%
\textit{, we have} 
\begin{eqnarray*}
&&\lim_{n\rightarrow \infty }\frac{1}{\widetilde{K}_{n}\left( x,x\right) ^{m}%
}R_{m}\left( x+\frac{\xi _{1}}{\widetilde{K}_{n}\left( x,x\right) },x+\frac{%
\xi _{2}}{\widetilde{K}_{n}\left( x,x\right) },...,x+\frac{\xi _{m}}{%
\widetilde{K}_{n}\left( x,x\right) }\right)  \\
&=&\det \left( \frac{\sin \pi \left( \xi _{i}-\xi _{j}\right) }{\pi \left(
\xi _{i}-\xi _{j}\right) }\right) _{i,j=1}^{m}.
\end{eqnarray*}%
\begin{equation}
\end{equation}%
\newline
\newline
\textbf{Corollary 1.3}\newline
\textit{Let }$r,s$\textit{\ be non-negative integers and }%
\begin{equation}
K_{n}^{\left( r,s\right) }\left( x,x\right) =\sum_{k=0}^{n-1}p_{k}^{\left(
r\right) }\left( x\right) p_{k}^{\left( r\right) }\left( s\right) .
\end{equation}%
\textit{Let  }%
\begin{equation}
\tau _{r,s}=\left\{ 
\begin{array}{rr}
0, & r+s\text{ odd} \\ 
\frac{\left( -1\right) ^{\left( r+s\right) /2}}{r+s+1}, & r+s\text{ even}%
\end{array}%
\right. .
\end{equation}%
\textit{Let }$I^{\prime }$\textit{\ be a closed subinterval of }$I^{0}.$ 
\textit{Then uniformly for} $x\in I^{\prime },$%
\begin{equation}
\lim_{n\rightarrow \infty }\frac{1}{n^{r+s+1}}K_{n}^{\left( r,s\right)
}\left( x,x\right) =\frac{1}{\pi w\left( x\right) \left( 1-x^{2}\right)
^{\left( r+s+1\right) /2}}\tau _{r,s}.
\end{equation}%
\newline
\textbf{Remarks}\newline
(a) We believe that the hypotheses above are the weakest imposed so far
guaranteeing universality for a fixed weight on $\left( -1,1\right) .$ Most
hypotheses imposed so far involve analyticity, for example in \cite%
{KuijlaarsVanlessen2002}.\newline
(b) The only reason for restricting $a,b$ to be real in (1.7), is that $%
\widetilde{K}_{n}\left( x+\frac{a}{\widetilde{K}_{n}\left( x,x\right) },x+%
\frac{b}{\widetilde{K}_{n}\left( x,x\right) }\right) $ involves the weight
evaluated at arguments involving $a$ and $b$. If we consider instead $%
K_{n}\left( x+\frac{a}{\widetilde{K}_{n}\left( x,x\right) },x+\frac{b}{%
\widetilde{K}_{n}\left( x,x\right) }\right) $, then the limits hold
uniformly for $a,b$ in compact subsets of the plane.

We also present $L_{p}$ results, assuming less about $w$:\newline
\newline
\textbf{Theorem 1.4}\newline
\textit{Let }$\mu $\textit{\ be a finite positive Borel measure on }$\left(
-1,1\right) $\textit{\ that is regular. Let }$p>0.$\textit{\ Let }$I$\textit{%
\ be a closed subinterval of }$\left( -1,1\right) $\textit{\ in which }$\mu $%
\textit{\ is absolutely continuous, and} $w$\textit{\ is bounded above and
below by positive constants, and moreover, }$w$\textit{\ is Riemann
integrable in }$I$\textit{. Then if }$I^{\prime }$ \textit{is a closed
subinterval of }$I^{0},$%
\begin{equation}
\lim_{n\rightarrow \infty }\int_{I^{\prime }}\left\vert \frac{\widetilde{K}%
_{n}\left( x+\frac{a}{\widetilde{K}_{n}\left( x,x\right) },x+\frac{b}{%
\widetilde{K}_{n}\left( x,x\right) }\right) }{\widetilde{K}_{n}\left(
x,x\right) }-\frac{\sin \pi \left( a-b\right) }{\pi \left( a-b\right) }%
\right\vert ^{p}dx=0,
\end{equation}%
\textit{\ uniformly for }$a,b$\textit{\ in compact subsets of the real line. 
\newline
}

The restriction of Riemann integrability of $w$ arises in showing $w\left( x+%
\frac{a}{\widetilde{K}_{n}\left( x,x\right) }\right) /w\left( x\right)
\rightarrow 1$ as $n\rightarrow \infty $, in a suitable sense. If we do not
assume $w$ is Riemann integrable in $I$, we can prove:\newline
\newline
\textbf{Theorem 1.5}\newline
\textit{Let }$\mu $\textit{\ be a finite positive Borel measure on }$\left(
-1,1\right) $\textit{\ that is regular. Let }$p>0.$\textit{\ Let }$I$\textit{%
\ be a closed subinterval of }$\left( -1,1\right) $\textit{\ in which }$\mu $%
\textit{\ is absolutely continuous, and} $w$\textit{\ is bounded above and
below by positive constants. Then if }$I^{\prime }$ \textit{is a closed
subinterval of }$I^{0}$, \textit{uniformly for }$a,b$ \textit{in compact
subsets of the plane,}%
\begin{equation}
\lim_{n\rightarrow \infty }\int_{I^{\prime }}\left\vert \frac{K_{n}\left( x+%
\frac{a}{\widetilde{K}_{n}\left( x,x\right) },x+\frac{b}{\widetilde{K}%
_{n}\left( x,x\right) }\right) }{K_{n}\left( x,x\right) }-\frac{\sin \pi
\left( a-b\right) }{\pi \left( a-b\right) }\right\vert ^{p}dx=0.
\end{equation}

When we assume only that $w$ is bounded below, and do not assume absolute
continuity of $\mu $, we can still prove an $L_{1}$ form of universality:%
\newline
\newline
\textbf{Theorem 1.6}\newline
\textit{Let }$\mu $\textit{\ be a finite positive Borel measure on }$\left(
-1,1\right) $\textit{\ that is regular. Let }$I$\textit{\ be a closed
subinterval of }$\left( -1,1\right) $\textit{\ in which }$w$\textit{\ is
bounded below by a positive constant. Then if }$I^{\prime }$ \textit{is a
closed subinterval of }$I^{0}$, \textit{uniformly for }$a,b$ \textit{in
compact subsets of the plane,}%
\begin{equation}
\lim_{n\rightarrow \infty }\int_{I^{\prime }}\left\vert \frac{1}{n}%
K_{n}\left( x+\frac{\pi a\sqrt{1-x^{2}}}{n},x+\frac{\pi b\sqrt{1-x^{2}}}{n}%
\right) -\frac{1}{\pi w\left( x\right) \sqrt{1-x^{2}}}\frac{\sin \pi \left(
a-b\right) }{\pi \left( a-b\right) }\right\vert dx=0.
\end{equation}

In the sequel\thinspace\ $C,C_{1},C_{2},...$ denote constants independent of 
$n,x,y,s,t$. The same symbol does not necessarily denote the same constant
in different occurences. We shall write $C=C\left( \alpha \right) $ or $%
C\neq C\left( \alpha \right) $ to respectively denote dependence on, or
independence of, the parameter $\alpha $. Given measures $\mu ^{\ast }$, $%
\mu ^{\#}$, we use $K_{n}^{\ast },K_{n}^{\#}$ and $p_{n}^{\ast },p_{n}^{\#}$
to denote respectively their reproducing kernels and orthonormal
polynomials. Similarly superscripts $\ast ,\#$ are used to distinguish other
quantities associated with them. The superscript $L$ denotes quantities
associated with the Legendre weight $1$ on $\left[ -1,1\right] $. For $x\in 
\mathbb{R}$ and $\delta >0$, we set%
\begin{equation*}
I\left( x,\delta \right) =\left[ x-\delta ,x+\delta \right] .
\end{equation*}%
Recall that the $n$th Christoffel function for a measure $\mu $ is 
\begin{equation*}
\lambda _{n}\left( x\right) =1/K_{n}\left( x,x\right) =\min_{\deg \left(
P\right) \leq n-1}\left( \int_{-1}^{1}P^{2}d\mu \right) /P^{2}\left(
x\right) .
\end{equation*}

The most important new idea in this paper is a localization principle for
universality. We use it repeatedly in various forms, but the following basic
inequality is typical. Suppose that $\mu ,\mu ^{\ast }$ are measures with $%
\mu \leq \mu ^{\ast }$ in $\left[ -1,1\right] $. Then for $x,y\in \left[ -1,1%
\right] ,$

\begin{eqnarray*}
&&\left\vert K_{n}\left( x,y\right) -K_{n}^{\ast }\left( x,y\right)
\right\vert /K_{n}\left( x,x\right) \\
&\leq &\left( \frac{K_{n}\left( y,y\right) }{K_{n}\left( x,x\right) }\right)
^{1/2}\left[ 1-\frac{K_{n}^{\ast }\left( x,x\right) }{K_{n}\left( x,x\right) 
}\right] ^{1/2} \\
&=&\left( \frac{\lambda _{n}\left( x\right) }{\lambda _{n}\left( y\right) }%
\right) ^{1/2}\left[ 1-\frac{\lambda _{n}\left( x\right) }{\lambda
_{n}^{\ast }\left( x\right) }\right] ^{1/2}.
\end{eqnarray*}%
Observe that on the right-hand side, we have only Christoffel functions, and
their asymptotics are very well understood.

The paper is organised as follows. In Section 2, we present some asymptotics
for Christoffel functions. In Section 3, we prove our localization
principle, including the above inequality. In Section 4, we approximate
locally the measure $\mu $ in Theorem 1.1 by a scaled Jacobi weight and then
prove Theorem 1.1. In Section 5, we prove the $L_{1}$ result Theorem 1.6,
and in Section 6, prove the $L_{p}$ results Theorem 1.4 and 1.5. In Section
7, we prove Corollaries 1.2 and 1.3.\newline
\newline
\textbf{Acknowledgement}\newline
This research was stimulated by the wonderful conference in honor of Percy
Deift's 60th birthday, held at the Courant Institute in June 2006. In the
present form, it was also inspired by a visit to Peter Sarnak at Princeton
University, and discussions with Eli Levin during our collaboration on \cite%
{LevinLubinsky2007}.

\section{Christoffel functions\newline
}

We use $\lambda _{n}^{L}$ to denote the $n$th Christoffel function for the
Legendre weight on $\left[ -1,1\right] $. The methods used to prove the
following result are very well known, but I could not find this theorem as
stated in the literature. The issue is that known asymptotics for
Christoffel functions do not include the increment $a/n$. We could use
existing results in \cite{Mateetal1991}, \cite{Nevai1979}, \cite{Nevai1986}, 
\cite{Totik2000} to treat the case where $x+a/n\in J$, and add a proof for
the case where this fails, but the amount of effort seems almost the same.%
\newline
\newline
\textbf{Theorem 2.1}\newline
\textit{Let }$\mu $\textit{\ be a regular measure on }$\left[ -1,1\right] .$%
\textit{\ Assume that }$\mu $\textit{\ is absolutely continuous in an open
interval containing }$J=\left[ c,d\right] $\textit{\ and in }$J$\textit{,} $%
w=\mu ^{\prime }$\textit{\ is positive and continuous.} \textit{Let }$A>0$. 
\textit{Then uniformly for }$a\in \left[ -A,A\right] ,$ \textit{and} $x\in J,
$ 
\begin{equation}
\lim_{n\rightarrow \infty }\lambda _{n}\left( x+\frac{a}{n}\right) /\lambda
_{n}^{L}\left( x+\frac{a}{n}\right) =w\left( x\right) .
\end{equation}%
\textit{Moreover, uniformly for }$n\geq n_{0}\left( A\right) ,x\in J$, 
\textit{and} $a\in \left[ -A,A\right] ,$%
\begin{equation}
\lambda _{n}\left( x+\frac{a}{n}\right) \sim \frac{1}{n}.
\end{equation}%
\textit{The constants implicit in }$\sim $\textit{\ do not depend on }$\rho $%
\textit{.} \newline
\textbf{Remarks}\newline
(a) The notation $\sim $ means that the ratio of the two Christoffel
functions is bounded above and below by positive constants independent of $n$
and $a$\newline
(b) We emphasize that we are assuming that $w$ is continuous in $\left[ c,d%
\right] $ when regarded as a function defined on $\left( -1,1\right) $.%
\newline
(c) Using asymptotics for $\lambda _{n}^{L}$, we can rewrite (2.1) as%
\begin{equation*}
\lim_{n\rightarrow \infty }n\lambda _{n}\left( x+\frac{a}{n}\right) =\pi 
\sqrt{1-x^{2}}w\left( x\right) .
\end{equation*}%
\newline
\textbf{Proof}\newline
Let $\varepsilon >0$ and choose $\delta >0$ such that $\mu $ is absolutely
continuous in 
\begin{equation*}
I=\left[ c-\delta ,d+\delta \right] \subset \left( -1,1\right) 
\end{equation*}%
and such that 
\begin{equation}
\left( 1+\varepsilon \right) ^{-1}\leq \frac{w\left( x\right) }{w\left(
y\right) }\leq 1+\varepsilon ,\text{ }x\in \left[ c-\delta ,d+\delta \right] 
\text{ with }\left\vert x-y\right\vert \leq \delta .
\end{equation}%
(Of course, this is possible because of uniform continuity and positivity of 
$w$). Let us fix $x_{0}\in J$, let%
\begin{equation*}
I\left( x_{0},\delta \right) =\left[ x_{0}-\delta ,x_{0}+\delta \right] 
\end{equation*}%
and define a measure $\mu ^{\ast }$ with 
\begin{equation*}
\mu ^{\ast }=\mu \text{ in }\left[ -1,1\right] \backslash I\left(
x_{0},\delta \right) 
\end{equation*}%
and in $I\left( x_{0},\delta \right) $, let $\mu ^{\ast }$ be absolutely
continuous, with absolutely continuous component $w^{\ast }$ satisfying 
\begin{equation}
w^{\ast }=w\left( x_{0}\right) \left( 1+\varepsilon \right) \text{ in }%
I\left( x_{0},\delta \right) .
\end{equation}%
Because of (2.3), $d\mu \leq d\mu ^{\ast }$ in $\left[ -1,1\right] ,$ so
that if $\lambda _{n}^{\ast }$ is the $n$th Christoffel function for $\mu
^{\ast }$, we have for all $x,$%
\begin{equation}
\lambda _{n}\left( x\right) \leq \lambda _{n}^{\ast }\left( x\right) .
\end{equation}%
We now find an upper bound for $\lambda _{n}^{\ast }\left( x\right) $\ for $%
x\in I\left( x_{0},\delta /2\right) $. There exists $r\in \left( 0,1\right) $
depending only on $\delta $ such that 
\begin{equation}
0\leq 1-\left( \frac{t-x}{2}\right) ^{2}\leq r\text{ for }x\in I\left(
x_{0},\delta /2\right) \text{ and }t\in \left[ -1,1\right] \backslash
I\left( x_{0},\delta \right) .
\end{equation}%
(In fact, we may take $r=1-\left( \frac{\delta }{4}\right) ^{2}$). Let $\eta
\in \left( 0,\frac{1}{2}\right) $ and choose $\sigma >1$ so close to $1$
that 
\begin{equation}
\sigma ^{1-\eta }<r^{-\eta /4}.
\end{equation}%
Let $m=m\left( n\right) =n-2\left[ \eta n/2\right] $. Fix $x\in I\left(
x_{0},\delta /2\right) $ and choose a polynomial $P_{m}$ of degree $\leq $ $%
m-1$ such that 
\begin{equation*}
\lambda _{m}^{L}\left( x\right) =\int_{-1}^{1}P_{m}^{2}\text{ and }%
P_{m}^{2}\left( x\right) =1.
\end{equation*}%
\newline
Thus $P_{m}$ is the minimizing polynomial in the Christoffel function for
the Legendre weight at $x$. Let 
\begin{equation*}
S_{n}\left( t\right) =P_{m}\left( t\right) \left( 1-\left( \frac{t-x}{2}%
\right) ^{2}\right) ^{\left[ \eta n/2\right] },
\end{equation*}%
a polynomial of degree $\leq $ $m-1+2\left[ \eta n/2\right] \leq n-1$ with $%
S_{n}\left( x\right) =1$. Then using (2.4) and (2.6), 
\begin{eqnarray*}
\lambda _{n}^{\ast }\left( x\right)  &\leq &\int_{-1}^{1}S_{n}^{2}d\mu
^{\ast } \\
&\leq &w\left( x_{0}\right) (1+\varepsilon )\int_{I\left( x_{0},\delta
\right) }P_{m}^{2}+\left\Vert P_{m}\right\Vert _{L_{\infty }\left( \left[
-1,1\right] \backslash I\left( x_{0},\delta \right) \right) }^{2}r^{\left[
\eta n/2\right] }\int_{\left[ -1,1\right] \backslash I\left( x_{0},\delta
\right) }d\mu ^{\ast } \\
&\leq &w\left( x_{0}\right) (1+\varepsilon )\lambda _{m}^{L}\left( x\right)
+\left\Vert P_{m}\right\Vert _{L_{\infty }\left[ -1,1\right] }^{2}r^{\left[
\eta n/2\right] }\int_{-1}^{1}d\mu ^{\ast }.
\end{eqnarray*}%
Now we use the key idea from \cite[Lemma 9, p. 450]{Mateetal1991}. For $%
m\geq m_{0}\left( \sigma \right) $, we have 
\begin{eqnarray*}
\left\Vert P_{m}\right\Vert _{L_{\infty }\left[ -1,1\right] }^{2} &\leq
&\sigma ^{m}\int_{-1}^{1}P_{m}^{2} \\
&=&\sigma ^{m}\lambda _{m}^{L}\left( x\right) .
\end{eqnarray*}%
(This holds more generally for any polynomial $P$ of degree $\leq m-1$, and
is a consequence of the regularity of the Legendre weight. Alternatively, we
could use classic bounds for the Christoffel functions for the Legendre
weight.) Then from (2.7), uniformly for $x\in I\left( x_{0},\delta /2\right) 
$,%
\begin{eqnarray*}
\lambda _{n}^{\ast }\left( x\right)  &\leq &w\left( x_{0}\right)
(1+\varepsilon )\lambda _{m}^{L}\left( x\right) \left\{ 1+C\left[ \sigma
^{1-\eta }r^{\eta /2}\right] ^{n}\right\}  \\
&\leq &w\left( x_{0}\right) (1+\varepsilon )\lambda _{m}^{L}\left( x\right)
\left\{ 1+o\left( 1\right) \right\} ,
\end{eqnarray*}%
so as $\lambda _{n}\leq \lambda _{n}^{\ast },$%
\begin{eqnarray}
&&\sup_{x\in I\left( x_{0},\delta /2\right) }\lambda _{n}\left( x\right)
/\lambda _{n}^{L}\left( x\right)   \notag \\
&\leq &w\left( x_{0}\right) (1+\varepsilon )\left\{ 1+o\left( 1\right)
\right\} \sup_{x\in I\left( x_{0},\delta \right) }\lambda _{m}^{L}\left(
x\right) /\lambda _{n}^{L}\left( x\right) .
\end{eqnarray}%
The $o\left( 1\right) $ terms is independent of $x_{0}$. Now for large
enough $n$, and some $C$ independent of $\eta ,m,n,x_{0},$%
\begin{equation}
\sup_{x\in \left[ -1,1\right] }\lambda _{m}^{L}\left( x\right) /\lambda
_{n}^{L}\left( x\right) \leq 1+C\eta .
\end{equation}%
Indeed if $\left\{ p_{k}^{L}\right\} $ denote the orthonormal Legendre
polynomials, they admit the bound \cite[p.170]{Nevai1979} 
\begin{equation*}
\left\vert p_{k}^{L}\left( x\right) \right\vert \leq C\left( 1-x^{2}+\frac{1%
}{k^{2}}\right) ^{-1/4},\text{ }x\in \left[ -1,1\right] .
\end{equation*}%
Then uniformly for $x\in \left[ -1,1\right] ,$ 
\begin{eqnarray*}
0 &\leq &1-\frac{\lambda _{n}^{L}\left( x\right) }{\lambda _{m}^{L}\left(
x\right) }=\lambda _{n}^{L}\left( x\right) \sum_{k=m}^{n-1}\left(
p_{k}^{L}\left( x\right) \right) ^{2} \\
&\leq &C\lambda _{n}^{L}\left( x\right) \left( n-m\right) \max_{\frac{n}{2}%
\leq k\leq n}\left( 1-x^{2}+\frac{1}{k^{2}}\right) ^{-1/2} \\
&\leq &C\eta n\lambda _{n}^{L}\left( x\right) \left( 1-x^{2}+\frac{1}{n^{2}}%
\right) ^{-1/2} \\
&\leq &C\eta ,
\end{eqnarray*}%
by classical bounds for Christoffel functions \cite[p. 108, Lemma 5]%
{Nevai1979}. So we have (2.9), and then (2.8) and (2.3) give for $n\geq
n_{0}=n_{0}\left( x_{0},\delta \right) ,$%
\begin{equation*}
\sup_{x\in I\left( x_{0},\delta /2\right) }\lambda _{n}\left( x\right)
/\left( \lambda _{n}^{L}\left( x\right) w\left( x\right) \right) \leq
(1+\varepsilon )^{2}(1+C\eta ).
\end{equation*}%
By covering $J$ with finitely many such intervals $I\left( x_{0},\delta
/2\right) $, we obtain for some maximal threshhold $n_{1}$, that for $n\geq
n_{1}=n_{1}\left( \varepsilon ,\delta ,J\right) ,$%
\begin{equation*}
\sup_{x\in \left[ c-\delta /2,d+\delta /2\right] }\lambda _{n}\left(
x\right) /\left( \lambda _{n}^{L}\left( x\right) w\left( x\right) \right)
\leq (1+\varepsilon )^{2}(1+C\eta ).
\end{equation*}%
It is is essential here that $C$ is independent of $\varepsilon ,\eta $. Now
let $A>0$ and $\left\vert a\right\vert \leq A$. There exists $%
n_{2}=n_{2}\left( A\right) $ such that for $n\geq n_{2}$ and all $\left\vert
a\right\vert \leq A$ and all $x\in J$, we have $x+\frac{a}{n}\in \left[
c-\delta /2,d+\delta /2\right] $. We deduce that 
\begin{equation*}
\limsup_{n\rightarrow \infty }\sup_{a\in \left[ -A,A\right] ,x\in J}\frac{%
\lambda _{n}\left( x+\frac{a}{n}\right) }{\lambda _{n}^{L}\left( x+\frac{a}{n%
}\right) w\left( x\right) }\leq (1+\varepsilon )^{2}\left( 1+C\eta \right) .
\end{equation*}%
As the left-hand side is independent of the parameters $\varepsilon ,\eta $,
we deduce that 
\begin{equation}
\limsup_{n\rightarrow \infty }\left( \sup_{a\in \left[ -A,A\right] ,x\in J}%
\frac{\lambda _{n}\left( x+\frac{a}{n}\right) }{\lambda _{n}^{L}\left( x+%
\frac{a}{n}\right) w\left( x\right) }\right) \leq 1.
\end{equation}%
In a similar way, we can establish the converse bound%
\begin{equation}
\limsup_{n\rightarrow \infty }\left( \sup_{a\in \left[ -A,A\right] ,x\in J}%
\frac{\lambda _{n}^{L}\left( x+\frac{a}{n}\right) w\left( x\right) }{\lambda
_{n}\left( x+\frac{a}{n}\right) }\right) \leq 1.
\end{equation}%
Indeed with $m,x$ and $\eta $ as above, let us choose a polynomial $P$ of
degree $\leq m-1$ such that%
\begin{equation*}
\lambda _{m}\left( x\right) =\int_{-1}^{1}P_{m}^{2}\left( t\right) d\mu
\left( t\right) \text{ and }P_{m}^{2}\left( x\right) =1.
\end{equation*}%
Then with $S_{n}$ as above, and proceeding as above,%
\begin{equation*}
\lambda _{n}^{L}\left( x\right) \leq \int_{-1}^{1}S_{n}^{2}
\end{equation*}%
\begin{eqnarray*}
&\leq &\left[ w\left( x_{0}\right) ^{-1}(1+\varepsilon )\right]
\int_{I\left( x_{0},\delta \right) }P_{m}^{2}d\mu +\left\Vert
P_{m}\right\Vert _{L_{\infty }\left( \left[ -1,1\right] \backslash I\left(
x_{0},\delta \right) \right) }^{2}r^{\left[ \eta n/2\right] }\int_{\left[
-1,1\right] \backslash I\left( x_{0},\delta \right) }1 \\
&\leq &\left[ w\left( x_{0}\right) ^{-1}(1+\varepsilon )\right] \lambda
_{m}\left( x\right) \left\{ 1+C\left[ \sigma ^{1-\eta }\left( 1-r\right)
^{\eta /2}\right] ^{n}\right\} ,
\end{eqnarray*}%
and so as above, 
\begin{eqnarray*}
&&\sup_{x\in I\left( x_{0},\delta /2\right) }\lambda _{m}^{L}\left( x\right)
/\lambda _{m}\left( x\right)  \\
&\leq &\left[ w\left( x_{0}\right) ^{-1}(1+\varepsilon )(1+o\left( 1\right) )%
\right] \sup_{x\in I\left( x_{0},\delta /2\right) }\lambda _{m}^{L}\left(
x\right) /\lambda _{n}^{L}\left( x\right)  \\
&\leq &\left[ w\left( x_{0}\right) ^{-1}(1+\varepsilon )\right] \left\{
1+o\left( 1\right) \right\} \left( 1+C\eta \right) .
\end{eqnarray*}%
Then (2.11) follows after a scale change $m\rightarrow n$ and using
monotonicity of $\lambda _{n}$ in \ $n$, much as above. Together (2.10) and
(2.11) give (2.1). Finally, (2.2) follows from standard bounds for the
Christoffel function for the Legendre weight. $\blacksquare $

\section{Localization}

\textbf{Theorem 3.1}\newline
\textit{Assume that }$\mu ,\mu ^{\ast }$\textit{\ are regular measures on }$%
\left[ -1,1\right] $ \textit{that are absolutely continuous in an open
interval containing }$J=\left[ c,d\right] $\textit{. Assume that }$w=\mu
^{\prime }$ \textit{is positive and continuous in }$J$\textit{\ and }%
\begin{equation*}
d\mu =d\mu ^{\ast }\text{ \textit{in} }J.
\end{equation*}%
\textit{Let }$A>0$\textit{. Then as }$n\rightarrow \infty ,$ 
\begin{equation}
\sup_{a,b\in \left[ -A,A\right] ,x\in J}\left\vert \left( K_{n}-K_{n}^{\ast
}\right) \left( x+\frac{a}{n},x+\frac{b}{n}\right) \right\vert /n=o\left(
1\right) .
\end{equation}%
\textbf{Proof\newline
}We initially assume that%
\begin{equation}
d\mu \leq d\mu ^{\ast }\text{ in }\left( -1,1\right) .
\end{equation}%
The idea is to estimate the $L_{2}$ norm of $K_{n}\left( x,t\right)
-K_{n}^{\ast }\left( x,t\right) $ over $\left[ -1,1\right] $, and then to
use Christoffel function estimates. Now 
\begin{eqnarray*}
&&\int_{-1}^{1}\left( K_{n}\left( x,t\right) -K_{n}^{\ast }\left( x,t\right)
\right) ^{2}d\mu \left( t\right)  \\
&=&\int_{-1}^{1}K_{n}^{2}\left( x,t\right) d\mu \left( t\right)
-2\int_{-1}^{1}K_{n}\left( x,t\right) K_{n}^{\ast }\left( x,t\right) d\mu
\left( t\right) +\int_{-1}^{1}K_{n}^{\ast 2}\left( x,t\right) d\mu \left(
t\right)  \\
&=&K_{n}\left( x,x\right) -2K_{n}^{\ast }\left( x,x\right)
+\int_{-1}^{1}K_{n}^{\ast 2}\left( x,t\right) d\mu \left( t\right) ,
\end{eqnarray*}%
by the reproducing kernel property. As $d\mu \leq d\mu ^{\ast },$ we also
have%
\begin{equation*}
\int_{-1}^{1}K_{n}^{\ast 2}\left( x,t\right) d\mu \left( t\right) \leq
\int_{-1}^{1}K_{n}^{\ast 2}\left( x,t\right) d\mu ^{\ast }\left( t\right)
=K_{n}^{\ast }\left( x,x\right) .
\end{equation*}%
So 
\begin{eqnarray}
&&\int_{-1}^{1}\left( K_{n}\left( x,t\right) -K_{n}^{\ast }\left( x,t\right)
\right) ^{2}d\mu \left( t\right)   \notag \\
&\leq &K_{n}\left( x,x\right) -K_{n}^{\ast }\left( x,x\right) .
\end{eqnarray}%
Next for any polynomial $P$ of degree $\leq n-1$, we have the Christoffel
function estimate 
\begin{equation}
\left\vert P\left( y\right) \right\vert \leq K_{n}\left( y,y\right)
^{1/2}\left( \int_{-1}^{1}P^{2}d\mu \right) ^{1/2}.
\end{equation}%
Applying this to $P\left( t\right) =K_{n}\left( x,t\right) -K_{n}^{\ast
}\left( x,t\right) $ and using (3.3) gives, for all $x,y\in \left[ -1,1%
\right] ,$ 
\begin{eqnarray*}
&&\left\vert K_{n}\left( x,y\right) -K_{n}^{\ast }\left( x,y\right)
\right\vert  \\
&\leq &K_{n}\left( y,y\right) ^{1/2}\left[ K_{n}\left( x,x\right)
-K_{n}^{\ast }\left( x,x\right) \right] ^{1/2}
\end{eqnarray*}%
so%
\begin{eqnarray}
&&\left\vert K_{n}\left( x,y\right) -K_{n}^{\ast }\left( x,y\right)
\right\vert /K_{n}\left( x,x\right)   \notag \\
&\leq &\left( \frac{K_{n}\left( y,y\right) }{K_{n}\left( x,x\right) }\right)
^{1/2}\left[ 1-\frac{K_{n}^{\ast }\left( x,x\right) }{K_{n}\left( x,x\right) 
}\right] ^{1/2}.
\end{eqnarray}%
Now we set $x=x_{0}+\frac{a}{n}$ and $y=x_{0}+\frac{b}{n}$, where $a,b\in %
\left[ -A,A\right] $ and $x_{0}\in J$. By Theorem 2.1, uniformly for such $x,
$ $\frac{K_{n}^{\ast }\left( x,x\right) }{K_{n}\left( x,x\right) }=1+o\left(
1\right) $, for they both have the same asymptotics as for the weight $w$ on
\ $\left[ -1,1\right] $. Moreover, uniformly for $a,b\in \left[ -A,A\right] ,
$ 
\begin{equation*}
K_{n}\left( x_{0}+\frac{b}{n},x_{0}+\frac{b}{n}\right) \sim K_{n}\left(
x_{0}+\frac{a}{n},x_{0}+\frac{a}{n}\right) \sim n,
\end{equation*}%
so 
\begin{equation*}
\sup_{a,b\in \left[ -A,A\right] ,x_{0}\in J}\left\vert \left(
K_{n}-K_{n}^{\ast }\right) \left( x_{0}+\frac{a}{n},x_{0}+\frac{b}{n}\right)
\right\vert /n=o\left( 1\right) .
\end{equation*}%
Now we drop the extra hypothesis (3.2). Define a measure $\nu $ by $\nu =\mu
=\mu ^{\ast }$ in $J;$ and in $\left[ -1,1\right] \backslash J$, let 
\begin{equation*}
d\nu \left( x\right) =\max \left\{ \left\vert x-c\right\vert \left\vert
x-d\right\vert ,w\left( x\right) ,w^{\ast }\left( x\right) \right\} dx+d\mu
_{s}\left( x\right) +d\mu _{s}^{\ast }\left( x\right) ,
\end{equation*}%
where $w,w^{\ast }$ and $\mu _{s},\mu _{s}^{\ast }$ are respectively the
absolutely continuous and singular components of $\mu ,\mu ^{\ast }$. Then $%
d\mu \leq d\nu $ and $d\mu ^{\ast }\leq d\nu $, and $\nu $ is regular as its
absolutely continuous component is positive in $\left( -1,1\right) $, and
hence lies in the even smaller class $\mathcal{M}.$ Moreover, $\nu $ is
absolutely continuous in an open interval containing $J,$ and $\nu ^{\prime
}=w$ in $J$ . The case above shows that the reproducing kernels for $\mu $
and $\mu ^{\ast }$ have the same asymptotics as that for $\nu $, in the
sense of (3.1), and hence the same asymptotics as each other. $\blacksquare $

\section{Smoothing}

In this section, we approximate $\mu $ of Theorem 1.1 by a scaled Legendre
Jacobi measure $\mu ^{\#}$ and then prove Theorem 1.1. Recall that $\tilde{K}%
_{n}$ is the normalized kernel, given by (1.5). Our smoothing result (which
may also be viewed as localization) is:\newline
\newline
\textbf{Theorem 4.1}\newline
\textit{Let }$\mu $\textit{\ be as in Theorem 1.1. Let }$A>0,\varepsilon \in
\left( 0,\frac{1}{2}\right) $\textit{\ and choose }$\delta >0$\textit{\ such
that (2.3) holds. Let }$x_{0}\in J.$ \textit{Then there exists }$C$ \textit{%
and} $n_{0}$\textit{\ such that for }$n\geq n_{0},$%
\begin{equation}
\sup_{a,b\in \left[ -A,A\right] ,x\in I\left( x_{o},\frac{\delta }{2}\right)
\cap J}\left\vert \left( \tilde{K}_{n}-K_{n}^{L}\right) \left( x+\frac{a}{n}%
,x+\frac{b}{n}\right) \right\vert /n\leq C\varepsilon ^{1/2},
\end{equation}%
\textit{where }$C$\textit{\ is independent of }$\varepsilon ,\delta ,n,x_{0}$%
.\newline
\textbf{Proof}\newline
Fix $x_{0}\in J$ and let $w^{\#}$ be the scaled Legendre weight%
\begin{equation*}
w^{\#}=w\left( x_{0}\right) \text{ in }\left( -1,1\right) .
\end{equation*}%
Note that 
\begin{equation}
K_{n}^{\#}\left( x,y\right) =\frac{1}{w\left( x_{0}\right) }K_{n}^{L}\left(
x,y\right) .
\end{equation}%
(Recall that the superscript indicates the Legendre weight on $\left[ -1,1%
\right] $). Because of our localization result Theorem 3.1, we may replace $%
d\mu $ by $w^{\ast }\left( x\right) dx$, where 
\begin{equation*}
w^{\ast }=w\text{ in }I\left( x_{0},\delta \right) 
\end{equation*}%
and 
\begin{equation*}
w^{\ast }=w\left( x_{0}\right) \text{ in }\left[ -1,1\right] \backslash
I\left( x_{0},\delta \right) ,
\end{equation*}%
without affecting the asymptotics for $K_{n}\left( x+\frac{a}{n},x+\frac{b}{n%
}\right) $ in the interval $I\left( x_{0},\frac{\delta }{2}\right) $. (Note
that $\varepsilon $ and $\delta $ play no role in Theorem 3.1). So in the
sequel, we assume that $w=w\left( x_{0}\right) =w^{\#}$ in $\left[ -1,1%
\right] \backslash I\left( x_{0},\delta \right) $, while not changing $w$ in 
$I$. Observe that (2.3) implies that 
\begin{equation}
\left( 1+\varepsilon \right) ^{-1}\leq \frac{w}{w^{\#}}\leq 1+\varepsilon 
\text{, in }\left[ -1,1\right] .
\end{equation}%
Then, much as in the previous section, 
\begin{eqnarray*}
&&\int_{-1}^{1}\left( K_{n}\left( x,t\right) -K_{n}^{\#}\left( x,t\right)
\right) ^{2}w^{\#}\left( t\right) dt \\
&=&\int_{-1}^{1}K_{n}^{2}\left( x,t\right) w^{\#}\left( t\right)
dt-2\int_{-1}^{1}K_{n}\left( x,t\right) K_{n}^{\#}\left( x,t\right)
w^{\#}\left( t\right) dt+\int_{-1}^{1}K_{n}^{\#2}\left( x,t\right)
w^{\#}\left( t\right) dt \\
&=&\int_{-1}^{1}K_{n}^{2}\left( x,t\right) w\left( t\right) dt+\int_{I\left(
x_{0},\delta \right) }K_{n}^{2}\left( x,t\right) \left( w^{\#}-w\right)
\left( t\right) dt-2K_{n}\left( x,x\right) +K_{n}^{\#}\left( x,x\right)  \\
&=&K_{n}^{\#}\left( x,x\right) -K_{n}\left( x,x\right) +\int_{I\left(
x_{0},\delta \right) }K_{n}^{2}\left( x,t\right) \left( w^{\#}-w\right)
\left( t\right) dt,
\end{eqnarray*}%
recall that $w=w^{\#}$ in $\left[ -1,1\right] \backslash I$. By (4.3), 
\begin{equation*}
\int_{I\left( x_{0},\delta \right) }K_{n}^{2}\left( x,t\right) \left(
w^{\#}-w\right) \left( t\right) dt\leq \varepsilon \int_{I\left(
x_{0},\delta \right) }K_{n}^{2}\left( x,t\right) w\left( t\right) dt\leq
\varepsilon K_{n}\left( x,x\right) .
\end{equation*}%
So%
\begin{equation}
\int_{-1}^{1}\left( K_{n}\left( x,t\right) -K_{n}^{\#}\left( x,t\right)
\right) ^{2}w^{\#}\left( t\right) dt\leq K_{n}^{\#}\left( x,x\right) -\left(
1-\varepsilon \right) K_{n}\left( x,x\right) .
\end{equation}%
Applying an obvious analogue of (3.4) to $P\left( t\right) =K_{n}\left(
x,t\right) -K_{n}^{\#}\left( x,t\right) $ and using (4.4) gives for $y\in %
\left[ -1,1\right] ,$ 
\begin{eqnarray*}
&&\left\vert K_{n}\left( x,y\right) -K_{n}^{\#}\left( x,y\right) \right\vert 
\\
&\leq &K_{n}^{\#}\left( y,y\right) ^{1/2}\left[ K_{n}^{\#}\left( x,x\right)
-\left( 1-\varepsilon \right) K_{n}\left( x,x\right) \right] ^{1/2}
\end{eqnarray*}%
so%
\begin{eqnarray*}
&&\left\vert K_{n}\left( x,y\right) -K_{n}^{\#}\left( x,y\right) \right\vert
/K_{n}^{\#}\left( x,x\right)  \\
&\leq &\left( \frac{K_{n}^{\#}\left( y,y\right) }{K_{n}^{\#}\left(
x,x\right) }\right) ^{1/2}\left[ 1-\left( 1-\varepsilon \right) \frac{%
K_{n}\left( x,x\right) }{K_{n}^{\#}\left( x,x\right) }\right] ^{1/2}.
\end{eqnarray*}%
In view of (4.3), we also have 
\begin{equation*}
\frac{K_{n}\left( x,x\right) }{K_{n}^{\#}\left( x,x\right) }=\frac{\lambda
_{n}^{\#}\left( x\right) }{\lambda _{n}\left( x\right) }\geq \frac{1}{%
1+\varepsilon },
\end{equation*}%
so for all $y\in \left[ -1,1\right] ,$%
\begin{eqnarray*}
&&\left\vert K_{n}\left( x,y\right) -K_{n}^{\#}\left( x,y\right) \right\vert
/K_{n}^{\#}\left( x,x\right)  \\
&\leq &\left( \frac{K_{n}^{\#}\left( y,y\right) }{K_{n}^{\#}\left(
x,x\right) }\right) ^{1/2}\left[ 1-\frac{1-\varepsilon }{1+\varepsilon }%
\right] ^{1/2} \\
&\leq &\sqrt{2\varepsilon }\left( \frac{K_{n}^{\#}\left( y,y\right) }{%
K_{n}^{\#}\left( x,x\right) }\right) ^{1/2} \\
&=&\sqrt{2\varepsilon }\left( \frac{K_{n}^{L}\left( y,y\right) }{%
K_{n}^{L}\left( x,x\right) }\right) ^{1/2}=\sqrt{2\varepsilon }\left( \frac{%
\lambda _{n}^{L}\left( x\right) }{\lambda _{n}^{L}\left( y\right) }\right)
^{1/2}.
\end{eqnarray*}%
Here we have used (4.2). Now we set $x=x_{1}+\frac{a}{n}$ and $y=x_{1}+\frac{%
b}{n}$, where $x_{1}\in I\left( x_{0},\frac{\delta }{2}\right) $ and $a,b\in %
\left[ -A,A\right] $. By classical estimates for Christoffel functions for
the Legendre weight (or even Theorem 2.1), uniformly for $a,b\in \left[ -A,A%
\right] ,$ and $x_{1}\in J,$ 
\begin{equation*}
\lambda _{n}^{L}\left( x_{1}+\frac{b}{n}\right) \sim \lambda _{n}^{L}\left(
x_{1}+\frac{a}{n}\right) \sim n^{-1},
\end{equation*}%
and also the constants implicit in $\sim $ are independent of $\varepsilon
,\delta $ and $x_{1}$ (this is crucial!). Thus for some $C$ and $n_{0}$
depending only on $A$ and $J$, we have for $n\geq n_{0},$%
\begin{equation*}
\sup_{a,b\in \left[ -A,A\right] ,x_{1}\in I\left( x_{0,}\frac{\delta }{2}%
\right) \cap J}\left\vert \left( K_{n}-K_{n}^{\#}\right) \left( x_{1}+\frac{a%
}{n},x_{1}+\frac{b}{n}\right) \right\vert /n\leq C\sqrt{\varepsilon }.
\end{equation*}%
Finally, recall (4.2), that%
\begin{equation*}
\left( 1+\varepsilon \right) ^{-1}\leq \frac{w\left( x_{1}\right) }{w\left(
x_{0}\right) }\leq 1+\varepsilon \text{ for all }x_{1}\in I\left(
x_{0},\delta \right) 
\end{equation*}%
and that $w$ is continuous at the endpoints of $I\left( x_{0,}\frac{\delta }{%
2}\right) \cap J$. $\blacksquare $\newline
\newline
\textbf{Proof of Theorem 1.1}\newline
Let $A,\varepsilon _{1}>0$. Choose $\varepsilon >0$ so small that the
right-hand side $C\varepsilon ^{1/2}$ of (4.1) is less than $\varepsilon _{1}
$. Choose $\delta >0$ such that (2.3) holds. Now divide $J$ into, say ~$M$
intervals $I\left( x_{j},\frac{\delta }{2}\right) $, $1\leq j\leq M$, each
of length $\delta $. For each $j$, there exists a threshhold $%
n_{0}=n_{0}\left( j\right) $ for which (4.1) holds for $n\geq n_{0}\left(
j\right) $ with $I\left( x_{0},\frac{\delta }{2}\right) $ replaced by $%
I\left( x_{j},\frac{\delta }{2}\right) $. Let $n_{1}$ denote the largest of
these. Then we obtain, for $n\geq n_{1}$, 
\begin{equation*}
\sup_{a,b\in \left[ -A,A\right] ,x_{0}\in J}\left\vert \left( \tilde{K}%
_{n}-K_{n}^{L}\right) \left( x+\frac{a}{n},x+\frac{b}{n}\right) \right\vert
/n\leq \varepsilon _{1}.
\end{equation*}%
It follows that 
\begin{equation}
\lim_{n\rightarrow \infty }\left( \sup_{a,b\in \left[ -A,A\right] ,x\in
J}\left\vert \left( \tilde{K}_{n}-K_{n}^{L}\right) \left( x+\frac{a}{n},x+%
\frac{b}{n}\right) \right\vert \right) =0.
\end{equation}%
Finally the universality limit for the Legendre weight (see for example \cite%
{KuijlaarsVanlessen2002}) gives as $n\rightarrow \infty $,%
\begin{eqnarray}
&&\frac{\pi \sqrt{1-x^{2}}}{n}K_{n}^{L}\left( x+\frac{u\pi \sqrt{1-x^{2}}}{n}%
,x+\frac{v\pi \sqrt{1-x^{2}}}{n}\right)   \notag \\
&\rightarrow &\frac{\sin \pi \left( u-v\right) }{\pi \left( u-v\right) },
\end{eqnarray}%
uniformly for $u,v$ in compact subsets of the real line, and $x$ in compact
subsets of $\left( -1,1\right) $. Setting 
\begin{equation*}
a=u\pi \sqrt{1-x^{2}}\text{ and }b=v\pi \sqrt{1-x^{2}}
\end{equation*}%
in (4.5), we obtain as $n\rightarrow \infty ,$ uniformly for $x\in J$ and $%
u,v$ in compact subsets of the real line, 
\begin{eqnarray}
&&\lim_{n\rightarrow \infty }\frac{\pi \sqrt{1-x^{2}}}{n}\tilde{K}_{n}\left(
x+\frac{u\pi \sqrt{1-x^{2}}}{n},x+\frac{v\pi \sqrt{1-x^{2}}}{n}\right)  
\notag \\
&=&\frac{\sin \pi \left( u-v\right) }{\pi \left( u-v\right) }.
\end{eqnarray}%
Since uniformly for $x\in J,$ by Theorem 2.1, 
\begin{eqnarray*}
\tilde{K}_{n}\left( x,x\right) ^{-1} &=&K_{n}^{L}\left( x,x\right)
^{-1}\left( 1+o\left( 1\right) \right)  \\
&=&\pi \sqrt{1-x^{2}}/n\left( 1+o\left( 1\right) \right) ,
\end{eqnarray*}%
we then also obtain the conclusion of Theorem 1.1. $\blacksquare $

For future use, we record also that 
\begin{equation}
\lim_{n\rightarrow \infty }\frac{1}{n}\tilde{K}_{n}\left( x+\frac{a}{n},x+%
\frac{b}{n}\right) =\frac{\sin \left( \left( a-b\right) /\sqrt{1-x^{2}}%
\right) }{\pi \left( a-b\right) }
\end{equation}%
uniformly for $x\in J$ and $a,b\in \left[ -A,A\right] $.

\section{Universality in $L_{1}$}

In this section, we prove Theorem 1.6. We assume that 
\begin{equation}
w\geq C_{0}\text{ in }I.
\end{equation}%
Let $\Delta >0$ and define a measure $\mu ^{\#}$ by 
\begin{equation*}
\mu ^{\#}=\mu \text{ in }\left[ -1,1\right] \backslash I
\end{equation*}%
and in $I$, we define $d\mu ^{\#}\left( x\right) =w^{\#}\left( x\right) dx$,
where%
\begin{equation}
w^{\#}\left( x\right) =\frac{1}{\Delta }\int_{x-\Delta }^{x+\Delta
}w=\int_{-1}^{1}w\left( x+s\Delta \right) ds.
\end{equation}%
\newline
\textbf{Lemma 5.1}\newline
\textit{Let }$I^{\prime }$\textit{\ be a closed subinterval of }$I^{0}$.%
\newline
\textit{(a) }$\mu ^{\#}$\textit{\ is continuous in }$I^{0}$\textit{\ and }$%
w^{\#}\geq \frac{1}{2}C_{0}$\textit{\ in }$I^{0}$\textit{.\newline
(b)} $\mu ^{\#}$ \textit{is regular on }$\left[ -1,1\right] $\textit{.}%
\newline
\textit{(c) There exists }$C_{1}>0$,\textit{\ independent of }$\Delta $,%
\textit{\ such that for} $n\geq 1,$ 
\begin{equation}
\sup_{t\in I^{\prime }}\frac{1}{n}K_{n}\left( t,t\right) \leq C_{1}\text{
and }\sup_{t\in I^{\prime }}\frac{1}{n}K_{n}^{\#}\left( t,t\right) \leq
C_{1}.
\end{equation}%
\textit{(d)} 
\begin{equation}
\lim_{n\rightarrow \infty }\frac{1}{n}\int_{I^{\prime }}\left\vert
K_{n}-K_{n}^{\#}\right\vert \left( t,t\right) dt=\frac{1}{\pi }%
\int_{I^{\prime }}\left\vert \frac{1}{w\left( t\right) }-\frac{1}{%
w^{\#}\left( t\right) }\right\vert \frac{dt}{\sqrt{1-t^{2}}}.
\end{equation}%
\textit{(e) For some }$C_{2}>0$\textit{\ independent of} $\Delta ,$%
\begin{eqnarray}
&&\int_{I^{\prime }}\frac{1}{\sqrt{1-t^{2}}}\left\vert \frac{1}{w\left(
t\right) }-\frac{1}{w^{\#}\left( t\right) }\right\vert dt  \notag \\
&\leq &C_{2}\sup_{\left\vert u\right\vert \leq \Delta }\int_{I}\left\vert
w\left( t+u\right) -w\left( t\right) \right\vert dt.
\end{eqnarray}%
\textbf{Proof}\newline
(a) is immediate.\newline
(b) This follows from Theorem 5.3.3 in \cite[p. 148]{StahlTotik1992}. As $%
\mu $ is regular, that theorem shows that the restriction of $\mu $ to $%
\left[ -1,1\right] \backslash I$ is regular. Hence the restriction of $\mu
^{\#}$ is trivially regular in $\left[ -1,1\right] \backslash I$. The
restriction of $\mu ^{\#}$ to $I$ is regular as its absolutely continuous
component $w^{\#}$\thinspace $>0$ there. Then Theorem 5.3.3 in \cite[p. 148]%
{StahlTotik1992} shows that $\mu ^{\#}$ is regular as a measure on all of $%
\left[ -1,1\right] .$\newline
(c) In view of (5.1), we have for $x\in I^{\prime },$ 
\begin{eqnarray*}
\lambda _{n}\left( x\right)  &\geq &C_{0}\inf_{\deg \left( P\right) \leq
n-1}\int_{I}P^{2}/P^{2}\left( x\right)  \\
&\geq &C_{0}C_{1}/n.
\end{eqnarray*}%
Here we are using classical bounds for the Legendre weight translated to the
interval $I$, and the constant $C_{1}$ depends only on the intervals $%
I^{\prime }$ and $I$. Then the first bound in (5.3) follows, and that for $%
\lambda _{n}^{\#}$ is similar. Since the lower bound on $\mu ^{\#}$ in $I$
is independent of $\Delta $, it follows that the constants we obtain in
(5.3) will also be independent of $\Delta .$\newline
(d) Since $\mu $ is regular, and $\mu ^{\prime }=w$ is bounded below by a
positive constant in $I$, we have a.e. in $I,$%
\begin{equation*}
\lim_{n\rightarrow \infty }K_{n}\left( x,x\right) /n=\frac{1}{\pi w\left(
x\right) \sqrt{1-x^{2}}}.
\end{equation*}%
See for example \cite[p. 449, Thm. 8]{Mateetal1991} or \cite[Theorem 1]%
{Totik2000}. A similar limit holds for $K_{n}^{\#}/n$. We also have the
uniform bound in (c). Then Lebesgue's Dominated Convergence Theorem gives
the result. \newline
(e) Recall that $I$ is a positive distance from $\pm 1$, while $w,w^{\#}$
are bounded below in $I$ by $C_{0}/2$. Then%
\begin{eqnarray*}
&&\int_{I^{\prime }}\frac{1}{\sqrt{1-t^{2}}}\left\vert \frac{1}{w\left(
t\right) }-\frac{1}{w^{\#}\left( t\right) }\right\vert dt \\
&\leq &C\int_{I^{\prime }}\left\vert w^{\#}\left( t\right) -w\left( t\right)
\right\vert dt \\
&\leq &C\int_{I^{\prime }}\int_{-1}^{1}\left\vert w\left( t+s\Delta \right)
-w\left( t\right) \right\vert ds\text{ }dt \\
&=&C\int_{-1}^{1}\int_{I^{\prime }}\left\vert w\left( t+s\Delta \right)
-w\left( t\right) \right\vert dt\text{ }ds \\
&\leq &C\sup_{\left\vert u\right\vert \leq \Delta }\int_{I^{\prime
}}\left\vert w\left( t+u\right) -w\left( t\right) \right\vert dt.
\end{eqnarray*}%
$\blacksquare $\newline
\newline
\textbf{Proof of Theorem 1.6}\newline
As per usual, 
\begin{eqnarray*}
&&\int_{-1}^{1}\left( K_{n}-K_{n}^{\#}\right) ^{2}\left( x,t\right) d\mu
^{\#}\left( t\right)  \\
&=&\int_{-1}^{1}K_{n}^{\#2}\left( x,t\right) d\mu ^{\#}\left( t\right)
-2\int_{-1}^{1}K_{n}^{\#}\left( x,t\right) K_{n}\left( x,t\right) d\mu
^{\#}\left( t\right) +\int_{-1}^{1}K_{n}^{2}\left( x,t\right) d\mu \left(
t\right)  \\
&&+\int_{I}K_{n}^{2}\left( x,t\right) d\left( \mu ^{\#}-\mu \right) \left(
t\right)  \\
&=&K_{n}^{\#}\left( x,x\right) -K_{n}\left( x,x\right)
+\int_{I}K_{n}^{2}\left( x,t\right) d\left( \mu ^{\#}-\mu \right) \left(
t\right)  \\
&\leq &K_{n}^{\#}\left( x,x\right) -K_{n}\left( x,x\right)
+\int_{I}K_{n}^{2}\left( x,t\right) \left( w^{\#}-w\right) \left( t\right) dt
\end{eqnarray*}%
recall that $\mu =\mu ^{\#}$ outside $I$ and that $\mu ^{\#}$ is absolutely
continuous in $I$. Then the Christoffel function estimate (3.4) gives for $%
y\in \left[ -1,1\right] ,$%
\begin{eqnarray*}
&&\left\vert K_{n}-K_{n}^{\#}\right\vert \left( x,y\right)  \\
&\leq &K_{n}^{\#}\left( y,y\right) ^{1/2}\left( K_{n}^{\#}\left( x,x\right)
-K_{n}\left( x,x\right) +\int_{I}K_{n}^{2}\left( x,t\right) \left(
w^{\#}-w\right) \left( t\right) dt\right) ^{1/2}.
\end{eqnarray*}%
\begin{equation}
\end{equation}%
We now replace $x$ by $x+\frac{a\pi \sqrt{1-x^{2}}}{n}$, $y$ by $x+\frac{%
a\pi \sqrt{1-x^{2}}}{n}$, integrate over $I^{\prime }$, and then use the
Cauchy-Schwarz inequality. We obtain 
\begin{eqnarray}
&&\int_{I^{\prime }}\left\vert K_{n}-K_{n}^{\#}\right\vert \left( x+\frac{%
a\pi \sqrt{1-x^{2}}}{n},x+\frac{b\pi \sqrt{1-x^{2}}}{n}\right) dx  \notag \\
&\leq &T_{1}^{1/2}T_{2}^{1/2},
\end{eqnarray}%
where 
\begin{equation*}
T_{1}=\int_{I^{\prime }}K_{n}^{\#}\left( x+\frac{b\pi \sqrt{1-x^{2}}}{n},x+%
\frac{b\pi \sqrt{1-x^{2}}}{n}\right) dx
\end{equation*}%
and%
\begin{eqnarray}
T_{2} &=&\int_{I^{\prime }}\left( K_{n}^{\#}-K_{n}\right) \left( x+\frac{%
a\pi \sqrt{1-x^{2}}}{n},x+\frac{a\pi \sqrt{1-x^{2}}}{n}\right) dx  \notag \\
&&+\int_{I^{\prime }}\left[ \int_{I}K_{n}^{2}\left( x+\frac{a\pi \sqrt{%
1-x^{2}}}{n},t\right) \left( w^{\#}-w\right) \left( t\right) dt\right] dx 
\notag \\
&=&:T_{21}+T_{22}.
\end{eqnarray}%
Now let $A>0$ and $a,b\in \left[ -A,A\right] $. Choose a subinterval $%
I^{\prime \prime }$ of $I^{0}$ such that $I^{\prime }\subset \left(
I^{\prime \prime }\right) ^{0}$. Observe that for some $n_{0}$ depending
only on $A$ and $I^{\prime },I^{\prime \prime }$, we have 
\begin{equation}
x+\frac{b\pi \sqrt{1-x^{2}}}{n}\in I^{\prime \prime }\text{ for }x\in
I^{\prime },b\in \left[ -A,A\right] ,n\geq n_{0}\text{.}
\end{equation}%
Then (c) of Lemma 5.1 shows that\ for $n\geq n_{0}$, 
\begin{equation}
T_{1}\leq C_{2}n,
\end{equation}%
where $C_{2}$ is independent of $n$ and $b\in \left[ -A,A\right] $. Next, we
make the substitution $s=x+\frac{a\pi \sqrt{1-x^{2}}}{n}$ in $T_{21}$.
Observe that 
\begin{equation*}
\frac{ds}{dx}=1-\frac{a\pi x}{n\sqrt{1-x^{2}}}\in \left[ \frac{1}{2},2\right]
,
\end{equation*}%
for $n\geq n_{1}$, where $n_{1}$ depends only on $A$ and $I$. We can also
assume that (5.9) holds for $n\geq n_{1}$. Hence for $n\geq \max \left\{
n_{0},n_{1}\right\} $ and all $a\in \left[ -A,A\right] ,$%
\begin{eqnarray*}
\left\vert T_{21}\right\vert  &\leq &\int_{I^{\prime }}\left\vert
K_{n}^{\#}-K_{n}\right\vert \left( x+\frac{a\pi \sqrt{1-x^{2}}}{n},x+\frac{%
a\pi \sqrt{1-x^{2}}}{n}\right) dx \\
&\leq &2\int_{I^{\prime \prime }}\left\vert K_{n}^{\#}-K_{n}\right\vert
\left( s,s\right) \text{ }ds
\end{eqnarray*}%
so using (d), (e) of the above lemma, 
\begin{equation*}
\limsup_{n\rightarrow \infty }\frac{1}{n}T_{21}\leq C\sup_{\left\vert
u\right\vert \leq \Delta }\int_{I^{\prime \prime }}\left\vert w\left(
t+u\right) -w\left( t\right) \right\vert dt,
\end{equation*}%
where $C$ does not depend on $\Delta $. Next, 
\begin{equation*}
\left\vert T_{22}\right\vert \leq \int_{I}\left\vert w-w^{\#}\right\vert
\left( t\right) \left[ \int_{I^{\prime }}K_{n}^{2}\left( x+\frac{a\pi \sqrt{%
1-x^{2}}}{n},t\right) dx\right] dt.
\end{equation*}%
Here for $n\geq \max \left\{ n_{0},n_{1}\right\} $, 
\begin{eqnarray*}
&&\int_{I^{\prime }}K_{n}^{2}\left( x+\frac{a\pi \sqrt{1-x^{2}}}{n},t\right)
dx \\
&\leq &\frac{1}{C_{0}}\int_{I^{\prime }}K_{n}^{2}\left( x+\frac{a\pi \sqrt{%
1-x^{2}}}{n},t\right) w\left( x+\frac{a\pi \sqrt{1-x^{2}}}{n}\right) dx \\
&\leq &\frac{2}{C_{0}}\int_{I^{\prime \prime }}K_{n}^{2}\left( s,t\right)
w\left( s\right) ds\leq \frac{2}{C_{0}}K_{n}\left( t,t\right) .
\end{eqnarray*}%
Then using\ (c), (e) of the previous lemma, we obtain%
\begin{eqnarray*}
\left\vert T_{22}\right\vert  &\leq &Cn\int_{I}\left\vert
w-w^{\#}\right\vert \left( t\right) dt \\
&\leq &Cn\sup_{\left\vert u\right\vert \leq \Delta }\int_{I^{\prime \prime
}}\left\vert w\left( t+u\right) -w\left( t\right) \right\vert dt.
\end{eqnarray*}%
Substituting all the above estimates in (5.7), we obtain%
\begin{eqnarray*}
&&\limsup_{n\rightarrow \infty }\frac{1}{n}\int_{I^{\prime }}\left\vert
K_{n}-K_{n}^{\#}\right\vert \left( x+\frac{a\pi \sqrt{1-x^{2}}}{n},x+\frac{%
b\pi \sqrt{1-x^{2}}}{n}\right) dx \\
&\leq &C\left( \sup_{\left\vert u\right\vert \leq \Delta }\int_{I^{\prime
\prime }}\left\vert w\left( t+u\right) -w\left( t\right) \right\vert
dt\right) ^{1/2},
\end{eqnarray*}%
uniformly for $a,b\in \left[ -A,A\right] $, where $C$ is independent of $%
\Delta $. Now as $\mu ^{\#}$ is regular, is absolutely continuous in $I$,
and $w^{\#}$ is continuous in $I^{0}$, Theorem 1.1 shows that 
\begin{eqnarray*}
&&\lim_{n\rightarrow \infty }\frac{1}{n}K_{n}^{\#}\left( x+\frac{a\pi \sqrt{%
1-x^{2}}}{n},x+\frac{b\pi \sqrt{1-x^{2}}}{n}\right)  \\
&=&\frac{\sin \pi \left( a-b\right) }{\pi \left( a-b\right) }\frac{1}{\pi 
\sqrt{1-x^{2}}w^{\#}\left( x\right) },
\end{eqnarray*}%
uniformly for $x\in I^{\prime }$ and $a,b\in \left[ -A,A\right] $. It
follows that 
\begin{eqnarray*}
&&\limsup_{n\rightarrow \infty }\int_{I^{\prime }}\left\vert \frac{1}{n}%
K_{n}\left( x+\frac{a\pi \sqrt{1-x^{2}}}{n},x+\frac{b\pi \sqrt{1-x^{2}}}{n}%
\right) -\frac{\sin \pi \left( a-b\right) }{\pi \left( a-b\right) }\frac{1}{%
\pi \sqrt{1-x^{2}}w\left( x\right) }\right\vert dx \\
&\leq &\left\vert \frac{\sin \pi \left( a-b\right) }{\pi \left( a-b\right) }%
\right\vert \int_{I^{\prime }}\frac{1}{\pi \sqrt{1-x^{2}}}\left\vert \frac{1%
}{w^{\#}\left( x\right) }-\frac{1}{w\left( x\right) }\right\vert dx \\
&&+C\left( \sup_{\left\vert u\right\vert \leq \Delta }\int_{I^{\prime \prime
}}\left\vert w\left( t+u\right) -w\left( t\right) \right\vert dt\right)
^{1/2},
\end{eqnarray*}%
uniformly for $a,b\in \left[ -A,A\right] $, where $C$ is independent of $%
\Delta $. Since the left-hand side is independent of $\Delta $, we may apply
(e) of the previous lemma, and then let $\Delta \rightarrow 0+$ to get the
result. Of course, as $w$ is integrable, we have as $\Delta \rightarrow 0+,$ 
\begin{equation*}
\sup_{\left\vert u\right\vert \leq \Delta }\int_{I^{\prime \prime
}}\left\vert w\left( t+u\right) -w\left( t\right) \right\vert dt\rightarrow
0.
\end{equation*}%
$\blacksquare $

\section{Universality in $L_{p}$}

The case $p=1$ of Theorem 1.5 is an immediate consequence of Theorem 1.6 and
the following lemma:\newline
\newline
\textbf{Lemma 6.1}\newline
\textit{Let }$A>0$\textit{\ and }$I^{\prime }$\textit{\ be a closed
subinterval of }$I^{0}$\textit{. As }$n\rightarrow \infty $, \textit{%
uniformly for }$a,b\in \left[ -A,A\right] $,%
\begin{eqnarray}
&&\frac{1}{n}\int_{I^{\prime }}\left\vert K_{n}\left( x+\frac{a\pi \sqrt{%
1-x^{2}}}{n},x+\frac{b\pi \sqrt{1-x^{2}}}{n}\right) -K_{n}\left( x+\frac{a}{%
\tilde{K}_{n}\left( x,x\right) },x+\frac{b}{\tilde{K}_{n}\left( x,x\right) }%
\right) \right\vert dx  \notag \\
&\rightarrow &0.
\end{eqnarray}%
\textbf{Proof}\newline
Choose a subinterval $I^{\prime \prime }$ of $I^{0}$ such that $I^{\prime
}\subset \left( I^{\prime \prime }\right) ^{0}$. Define $r_{n}\left(
x\right) $ by 
\begin{equation*}
\frac{1}{\tilde{K}_{n}\left( x,x\right) }=\frac{a\pi \sqrt{1-x^{2}}}{n}%
r_{n}\left( x\right) .
\end{equation*}%
Then the integrand in (6.1) may be written as 
\begin{eqnarray*}
&&\left\vert 
\begin{array}{c}
K_{n}\left( x+\frac{a\pi \sqrt{1-x^{2}}}{n},x+\frac{b\pi \sqrt{1-x^{2}}}{n}%
\right)  \\ 
-K_{n}\left( x+\frac{a\pi \sqrt{1-x^{2}}}{n}r_{n}\left( x\right) ,x+\frac{%
b\pi \sqrt{1-x^{2}}}{n}r_{n}\left( x\right) \right) 
\end{array}%
\right\vert  \\
&\leq &\left\vert \frac{\partial }{\partial s}K_{n}\left( s,x+\frac{b\pi 
\sqrt{1-x^{2}}}{n}\right) \right\vert _{|s=\xi }\frac{\left\vert
a\right\vert \pi \sqrt{1-x^{2}}}{n}\left\vert 1-r_{n}\left( x\right)
\right\vert  \\
&&+\left\vert \frac{\partial }{\partial t}K_{n}\left( x+\frac{a\pi \sqrt{%
1-x^{2}}}{n}r_{n}\left( x\right) ,t\right) \right\vert _{|t=\zeta }\frac{%
\left\vert b\right\vert \pi \sqrt{1-x^{2}}}{n}\left\vert 1-r_{n}\left(
x\right) \right\vert 
\end{eqnarray*}%
where $\xi $ lies between $x+\frac{a\pi \sqrt{1-x^{2}}}{n}$ and $x+\frac{%
a\pi \sqrt{1-x^{2}}}{n}r_{n}\left( x\right) $, with a similar restriction on 
$\zeta $. Now by Lemma 5.1(c) and Cauchy-Schwarz, 
\begin{equation*}
\sup_{s,t\in I}\left\vert K_{n}\left( s,t\right) \right\vert \leq Cn.
\end{equation*}%
By Bernstein's inequality \cite[p. 98, Corollary 1.2]{DeVoreLorentz1993}, 
\begin{equation*}
\sup_{s\in I^{\prime \prime },t\in I}\left\vert \frac{\partial }{\partial s}%
K_{n}\left( s,t\right) \right\vert \leq C_{1}Cn
\end{equation*}%
with a similar bound for $\frac{\partial }{\partial t}K_{n}$. Here $C_{1}$
depends only on $I$ and $I^{\prime \prime }$. Then for some $C_{2}$
independent of $a,b,n,x,$%
\begin{eqnarray*}
&&\frac{1}{n}\left\vert 
\begin{array}{c}
K_{n}\left( x+\frac{a\pi \sqrt{1-x^{2}}}{n},x+\frac{b\pi \sqrt{1-x^{2}}}{n}%
\right)  \\ 
-K_{n}\left( x+\frac{a\pi \sqrt{1-x^{2}}}{n}r_{n}\left( x\right) ,x+\frac{%
b\pi \sqrt{1-x^{2}}}{n}r_{n}\left( x\right) \right) 
\end{array}%
\right\vert  \\
&\leq &C\left\vert 1-r_{n}\left( x\right) \right\vert .
\end{eqnarray*}%
\newline
Hence the integral in the left-hand side of (6.1) is bounded above by 
\begin{equation*}
C\int_{I^{\prime }}\left\vert 1-r_{n}\left( x\right) \right\vert dx.
\end{equation*}%
Of course $C$ is independent of $n$. Next \cite[p. 449, Thm. 8]{Mateetal1991}%
, 
\begin{equation}
r_{n}\left( x\right) =\frac{n}{K_{n}\left( x,x\right) w\left( x\right) \pi 
\sqrt{1-x^{2}}}\rightarrow 1\text{ a.e. in }I
\end{equation}%
by Theorem 2.1. We shall shortly show that \ 
\begin{equation}
r_{n}\left( x\right) \leq C\text{ for }x\in I^{\prime }\text{ and }n\geq
n_{0}\text{.}
\end{equation}%
Then Lebesgue's Dominated Convergence Theorems shows that 
\begin{equation*}
\lim_{n\rightarrow \infty }\int_{I^{\prime }}\left\vert 1-r_{n}\left(
x\right) \right\vert dx=0.
\end{equation*}%
To prove (6.3), choose $M>0$ such that $w\leq M$ in $I$. Define a measure $%
\mu ^{\ast }$ by 
\begin{eqnarray*}
d\mu  &=&d\mu ^{\ast }\text{ in }\left[ -1,1\right] \backslash I; \\
d\mu ^{\ast }\left( x\right)  &=&Mdx\text{ in }I.
\end{eqnarray*}%
Then $d\mu \leq d\mu ^{\ast }$ in $\left[ -1,1\right] $ so $\lambda _{n}\leq
\lambda _{n}^{\ast }$ in $\left[ -1,1\right] $. As the absolutely continuous
component of $\mu ^{\ast }$ is positive and continuous in $I$, Theorem 2.1
shows that for some $C>0$, 
\begin{equation*}
\lambda _{n}^{\ast }\left( x\right) \leq \frac{C}{n}\text{ for }x\in
I^{\prime }\text{ and }n\geq 1,
\end{equation*}%
and then 
\begin{equation}
\frac{n}{K_{n}\left( x,x\right) }=n\lambda _{n}\left( x\right) \leq C\text{
for }x\in I^{\prime }\text{ and }n\geq 1.
\end{equation}%
The definition (6.2) of $r_{n}$, the fact that $w$ is bounded below in $I$,
and this last inequality, give (6.3). $\blacksquare $\newline
\newline
\textbf{Proof of Theorem 1.5}\newline
As $w$ is bounded above and below in $I$, the lemma and Theorem 1.6 give
uniformly for $a,b\in \left[ -A,A\right] ,$%
\begin{equation*}
\lim_{n\rightarrow \infty }\int_{I^{\prime }}\left\vert K_{n}\left( x+\frac{a%
}{\tilde{K}_{n}\left( x,x\right) },x+\frac{b}{\tilde{K}_{n}\left( x,x\right) 
}\right) \frac{w\left( x\right) \pi \sqrt{1-x^{2}}}{n}-\frac{\sin \pi \left(
a-b\right) }{\pi \left( a-b\right) }\right\vert dx=0.
\end{equation*}%
Now a.e. in $I,$ 
\begin{equation*}
\frac{1}{K_{n}\left( x,x\right) }=\frac{w\left( x\right) \pi \sqrt{1-x^{2}}}{%
n}\left( 1+o\left( 1\right) \right) .
\end{equation*}%
Moreover, by (6.4), Lemma 5.1(c), and Cauchy-Schwarz, both $\frac{1}{n}%
K_{n}\left( x+\frac{a}{\tilde{K}_{n}\left( x,x\right) },x+\frac{b}{\tilde{K}%
_{n}\left( x,x\right) }\right) $ and $K_{n}\left( x+\frac{a}{\tilde{K}%
_{n}\left( x,x\right) },x+\frac{b}{\tilde{K}_{n}\left( x,x\right) }\right)
/K_{n}\left( x,x\right) $ are bounded above uniformly for $a,b\in \left[ -A,A%
\right] ,$ $x\in I^{\prime }$, and $n\geq n_{0}$. We deduce that 
\begin{equation*}
\lim_{n\rightarrow \infty }\int_{I^{\prime }}\left\vert K_{n}\left( x+\frac{a%
}{\tilde{K}_{n}\left( x,x\right) },x+\frac{b}{\tilde{K}_{n}\left( x,x\right) 
}\right) /K_{n}\left( x,x\right) -\frac{\sin \pi \left( a-b\right) }{\pi
\left( a-b\right) }\right\vert dx=0.
\end{equation*}%
Finally, as we have just noted, the integrand in the last integral is
bounded above uniformly for $a,b\in \left[ -A,A\right] ,$ $x\in I^{\prime }$%
, and $n\geq n_{0}$, so we may replace the first power by the $p$th power,
for any $p>1$. For $p<1$, we can use H\"{o}lder's inequality. $\blacksquare $%
\newline
\newline
In proving Theorem 1.4, our last step is to replace $\frac{K_{n}\left( x+%
\frac{a}{\widetilde{K}_{n}\left( x,x\right) },x+\frac{b}{\widetilde{K}%
_{n}\left( x,x\right) }\right) }{K_{n}\left( x,x\right) }$ by $\frac{%
\widetilde{K}_{n}\left( x+\frac{a}{\widetilde{K}_{n}\left( x,x\right) },x+%
\frac{b}{\widetilde{K}_{n}\left( x,x\right) }\right) }{\widetilde{K}%
_{n}\left( x,x\right) }$. This is more difficult than one might expect - it
is only here that we need Riemann integrability of $w$ in $I$. For general
Lebesgue measurable $w$, it seems difficult to deal with the factor $%
\widetilde{K}_{n}\left( x,x\right) =w\left( x\right) K_{n}\left( x,x\right) $
below.\newline
\newline
\textbf{Lemma 6.2} \newline
\textit{Assume that }$w$\textit{\ is Riemann integrable and bounded below by
a positive constant in }$I$\textit{. Let }$I^{\prime }$\textit{\ be a
compact subinterval of }$I$\textit{. Let }$p,A>0$\textit{. Then uniformly
for }$a,b\in \left[ -A,A\right] $\textit{, we have} 
\begin{equation*}
\lim_{n\rightarrow \infty }\int_{I^{\prime }}\left\vert \sqrt{w\left( x+%
\frac{a}{\widetilde{K}_{n}\left( x,x\right) }\right) w\left( x+\frac{b}{%
\widetilde{K}_{n}\left( x,x\right) }\right) }/w\left( x\right) -1\right\vert
^{p}dx=0.
\end{equation*}%
\textbf{Proof}\newline
Let $a,b\in \left[ -A,A\right] $. From (6.4), for a suitable integer $n_{0}$
and some $L>0$, we have 
\begin{equation*}
\left\vert \frac{a}{\widetilde{K}_{n}\left( x,x\right) }\right\vert \leq 
\frac{L}{n}\text{ and }\left\vert \frac{b}{\widetilde{K}_{n}\left(
x,x\right) }\right\vert \leq \frac{L}{n},
\end{equation*}%
uniformly for $x\in I^{\prime }$, $a,b\in \left[ -A,A\right] $, and $n\geq
n_{0}$. Next, as $w$ is Riemann integrable in $I$, it is continuous a.e. in $%
I$ \cite[p. 23]{RieszNagy1990}. For $x\in I$ and $n\geq 1,$ let 
\begin{equation*}
\Omega _{n}\left( x\right) =\sup \left\{ \left\vert w\left( x+s\right)
-w\left( x\right) \right\vert :\left\vert s\right\vert \leq \frac{L}{n}%
\right\} .
\end{equation*}%
Note that for $x\in I^{\prime },n\geq n_{0}$ and $a,b\in \left[ -A,A\right] $%
, 
\begin{equation*}
\left\vert w\left( x+\frac{a}{\widetilde{K}_{n}\left( x,x\right) }\right)
-w\left( x\right) \right\vert \leq \Omega _{n}\left( x\right) .
\end{equation*}%
We have at every point of continuity of $w$ and in particular for a.e. $x\in
I$, 
\begin{equation*}
\lim_{n\rightarrow \infty }\Omega _{n}\left( x\right) =0.
\end{equation*}%
Moreover, as $w$ is Riemann integrable, $\Omega _{n}$ is bounded above in $I$%
, uniformly in $n$. Then Lebesgue's Dominated Convergence Theorem gives
uniformly for $a\in \left[ -A,A\right] ,$%
\begin{eqnarray*}
&&\int_{I^{\prime }}\left\vert w\left( x+\frac{a}{\widetilde{K}_{n}\left(
x,x\right) }\right) -w\left( x\right) \right\vert ^{p}dx \\
&\leq &\int_{I^{\prime }}\Omega _{n}\left( x\right) ^{p}dx\rightarrow 0\text{%
, }n\rightarrow \infty .
\end{eqnarray*}%
This, the fact that $w$ is bounded above and below, and some elementary
manipulations, give the result. $\blacksquare $\newline
\newline
\textbf{Proof of Theorem 1.4}\newline
Since $\frac{K_{n}\left( x+\frac{a}{\widetilde{K}_{n}\left( x,x\right) },x+%
\frac{b}{\widetilde{K}_{n}\left( x,x\right) }\right) }{K_{n}\left(
x,x\right) }$ is bounded uniformly in $n,x,a,b$ (over the relevant ranges)
and 
\begin{eqnarray*}
&&\frac{\widetilde{K}_{n}\left( x+\frac{a}{\widetilde{K}_{n}\left(
x,x\right) },x+\frac{b}{\widetilde{K}_{n}\left( x,x\right) }\right) }{%
\widetilde{K}_{n}\left( x,x\right) }/\frac{K_{n}\left( x+\frac{a}{\widetilde{%
K}_{n}\left( x,x\right) },x+\frac{b}{\widetilde{K}_{n}\left( x,x\right) }%
\right) }{K_{n}\left( x,x\right) } \\
&=&\sqrt{w\left( x+\frac{a}{\widetilde{K}_{n}\left( x,x\right) }\right)
w\left( x+\frac{b}{\widetilde{K}_{n}\left( x,x\right) }\right) }/w\left(
x\right) ,
\end{eqnarray*}%
this follows directly from the lemma above and Theorem 1.5. $\blacksquare $

\section{Proof of Corollaries 1.2 and 1.3}

\textbf{Proof of Corollary 1.2}\newline
This follows directly by substituting (1.6) into the determinant defining $%
R_{m}$. $\blacksquare $\newline
\newline
In proving Corollary 1.3, we need \newline
\newline
\textbf{Lemma 7.1}\newline
\textit{Let }$w\geq C$ \textit{in} $I$\textit{\ \ and }$I^{\prime
},I^{\prime \prime }$\textit{\ be closed subintervals of }$I^{0}$\textit{\
such that }$I^{\prime }$\textit{\ is contained in the interior of} $%
I^{\prime \prime }$. \textit{Let }$A>0$. \textit{There exists }$C_{2}$%
\textit{\ such that for }$n\geq 1,x\in I^{\prime },$\textit{\ and all }$%
\alpha ,\beta \in \mathbb{C}$\textit{\ with }$\left\vert \alpha \right\vert
,\left\vert \beta \right\vert \leq A$\textit{, }%
\begin{equation}
\left\vert \frac{1}{n}K_{n}\left( x+\frac{\alpha }{n},x+\frac{\beta }{n}%
\right) \right\vert \leq C_{2}.
\end{equation}%
\textit{\newline
}\textbf{Proof}\newline
Recall that $\frac{1}{n}K_{n}\left( x,x\right) $ is uniformly bounded above
for $x\in I^{\prime }$ by Lemma 5.1(c). Applying Cauchy-Schwarz, we obtain
for $x,y\in I^{\prime \prime },$%
\begin{equation}
\frac{1}{n}\left\vert K_{n}\left( x,y\right) \right\vert \leq \sqrt{\frac{1}{%
n}K_{n}\left( x,x\right) }\sqrt{\frac{1}{n}K_{n}\left( y,y\right) }\leq
C_{1}.
\end{equation}%
Next we note Bernstein's growth lemma for polynomials in the plane \cite[%
Theorem 2.2, p. 101]{DeVoreLorentz1993}: if $P$ is a polynomial of degree $%
\leq n$, we have for $z\notin \left[ -1,1\right] $,%
\begin{equation*}
\left\vert P\left( z\right) \right\vert \leq \left\vert z+\sqrt{z^{2}-1}%
\right\vert ^{n}\left\Vert P\right\Vert _{L_{\infty }\left[ -1,1\right] }.
\end{equation*}%
From this we deduce that given $L>0$, and $0<\delta <1$, there exists $%
C_{2}\neq C_{2}\left( n,P,z\right) $ such that for $\left\vert \func{Re}%
\left( z\right) \right\vert \leq \delta $, and $\left\vert \func{Im}%
z\right\vert \leq \frac{L}{n}$ 
\begin{equation*}
\left\vert P\left( z\right) \right\vert \leq C_{2}\left\Vert P\right\Vert
_{L_{\infty }\left[ -1,1\right] }.
\end{equation*}%
Mapping this to $I$ by a linear transformation, we deduce that for $\func{Re}%
z\in I^{\prime }$ and $\left\vert \func{Im}z\right\vert \leq \frac{L}{n},$ 
\begin{equation*}
\left\vert P\left( z\right) \right\vert \leq C_{3}\left\Vert P\right\Vert
_{L_{\infty }\left( I^{\prime \prime }\right) }
\end{equation*}%
where $C_{3}\neq C_{3}\left( n,P,z\right) $. We now apply this to $\frac{1}{n%
}K_{n}\left( x,y\right) $, separately in each variable, obtaining the stated
result. $\blacksquare $\newline
\newline
\textbf{Proof of Corollary 1.3}\newline
Note first from the lemma, $\left\{ \frac{1}{n}K_{n}\left( x+\frac{\alpha }{n%
},x+\frac{\beta }{n}\right) \right\} _{n=1}^{\infty }$ is analytic in $%
\alpha ,\beta $ and uniformly bounded for $\alpha ,\beta $ in compact
subsets of the plane. Moreover, from (4.8), and continuity of $w$,%
\begin{equation*}
\lim_{n\rightarrow \infty }\frac{1}{n}w\left( x\right) K_{n}\left( x+\frac{%
\alpha }{n},x+\frac{\beta }{n}\right) =\frac{\sin \left( \left( \alpha
-\beta \right) /\sqrt{1-x^{2}}\right) }{\pi \left( \alpha -\beta \right) }
\end{equation*}%
uniformly for $x\in I^{\prime }$ and $\alpha ,\beta $ in compact subsets of $%
I^{\prime }$. By convergence continuation theorems, this last limit then
holds uniformly for $\alpha ,\beta $ in compact subsets of the plane. Next,
expanding $p_{k}\left( x+\frac{\alpha }{n}\right) $ and $p_{k}\left( x+\frac{%
\beta }{n}\right) $ in Taylor series about $x,$ 
\begin{eqnarray*}
\frac{1}{n}K_{n}\left( x+\frac{\alpha }{n},x+\frac{\beta }{n}\right)  &=&%
\frac{1}{n}\sum_{k=0}^{n-1}p_{k}\left( x+\frac{\alpha }{n}\right)
p_{k}\left( x+\frac{\beta }{n}\right)  \\
&=&\frac{1}{n}\sum_{r,s=0}^{\infty }\frac{\left( \frac{\alpha }{n}\right)
^{r}}{r!}\frac{\left( \frac{\beta }{n}\right) ^{s}}{s!}%
\sum_{k=0}^{n-1}p_{k}^{\left( r\right) }\left( x\right) p_{k}^{\left(
s\right) }\left( x\right)  \\
&=&\sum_{r,s=0}^{\infty }\frac{\alpha ^{r}}{r!}\frac{\beta ^{s}}{s!}\frac{1}{%
n^{r+s+1}}K_{n}^{\left( r,s\right) }\left( x,x\right) ,
\end{eqnarray*}%
with the notation (1.9). Since the series terminates, the interchanges are
valid. By using the Maclaurin series of $\sin $ and the binomial theorem, we
see that 
\begin{equation*}
\frac{\sin \left( \alpha -\beta \right) }{\alpha -\beta }=\sum_{r,s=0}^{%
\infty }\frac{\alpha ^{r}}{r!}\frac{\beta ^{s}}{s!}\tau _{r,s},
\end{equation*}%
where $\tau _{r,s}$ is given by (1.10). Since uniformly convergent sequences
of analytic functions have Taylor series that also converge, we see that for 
$x\in I$, and each $r,s\geq 0,$ 
\begin{equation*}
\lim_{n\rightarrow \infty }\frac{1}{n^{r+s+1}}w\left( x\right) K_{n}^{\left(
r,s\right) }\left( x,x\right) =\frac{\tau _{r,s}}{\pi }\left( 1-x^{2}\right)
^{-\left( r+s\right) /2}.
\end{equation*}%
This establishes the limit (1.11), but we must still prove uniformity. Let $%
A,\varepsilon >0$. By the uniform convergence in Theorem 1.1, there exists $%
n_{0}$ such that for $n\geq n_{0},$%
\begin{eqnarray}
&&\left\vert 
\begin{array}{c}
\frac{w\left( x\right) \sqrt{1-x^{2}}}{n}K_{n}\left( x+\frac{a\pi \sqrt{%
1-x^{2}}}{n},x+\frac{b\pi \sqrt{1-x^{2}}}{n}\right)  \\ 
-\frac{w\left( y\right) \sqrt{1-y^{2}}}{n}K_{n}\left( y+\frac{a\pi \sqrt{%
1-y^{2}}}{n},y+\frac{b\pi \sqrt{1-y^{2}}}{n}\right) 
\end{array}%
\right\vert   \notag \\
&\leq &\varepsilon ,
\end{eqnarray}%
uniformly for $x,y\in I,a,b\in \left[ -A,A\right] $ and $n\geq n_{0}$. Using
Bernstein's growth inequality as in the lemma above, applied to the
polynomial in the left-hand side of (7.3), we obtain that this inequality
persists for complex $\alpha ,\beta $ with $\left\vert \alpha \right\vert
,\left\vert \beta \right\vert \leq A$, except that we must replace $%
\varepsilon $ by $C\varepsilon $, where $C$ depends only on $A$, not on $%
n,x,a,b,\varepsilon $. We can now use Cauchy's inequalities to bound the
Taylor series coefficients of the double series in $a,b$ implicit in the
left-hand side in (7.3). This leads to bounds on%
\begin{equation*}
\left\vert \frac{1}{n^{r+s+1}}w\left( x\right) K_{n}^{\left( r,s\right)
}\left( x,x\right) -\frac{1}{n^{r+s+1}}w\left( y\right) K_{n}^{\left(
r,s\right) }\left( y,y\right) \right\vert 
\end{equation*}%
that are uniform in $x,y$. $\blacksquare $

\end{document}